\documentclass[reqno]{amsart} 

\usepackage{amsmath,amsfonts,amssymb,amsthm}

\usepackage[T1]{fontenc}

\usepackage{etoolbox}
\patchcmd{\section}{\scshape}{\scshape\bfseries}{}{}
\makeatletter
\renewcommand{\@secnumfont}{\scshape\bfseries}
\makeatother

\usepackage{parskip}

\usepackage[bookmarks=true,bookmarksopen=true]{hyperref}

\theoremstyle{theorem}
\newtheorem{theorem}{Theorem}    

\newtheorem{conjecture}[theorem]{Conjecture}
\newtheorem*{question}{Open Question}
\newtheorem{corollary}[theorem]{Corollary}
\newtheorem{lemma}[theorem]{Lemma}
\newtheorem{proposition}[theorem]{Proposition}

\numberwithin{theorem}{section}

\theoremstyle{plain}
\newtheorem{Atheorem}{Theorem}
\newtheorem{Acorollary}[Atheorem]{Corollary}
\newtheorem{Aproposition}[Atheorem]{Proposition}

\theoremstyle{definition}
\newtheorem*{definition}{Definition}
\newtheorem{example}[theorem]{Example}

\theoremstyle{remark}

\numberwithin{remark}{section}

\numberwithin{equation}{section}

\numberwithin{table}{section}

\usepackage{tikz}
\usetikzlibrary{matrix}

\usepackage{todonotes}

\usepackage{enumerate}

\usepackage{float}

\makeatletter
\def\blfootnote{\gdef\@thefnmark{}\@footnotetext}
\makeatother

\DeclareMathOperator{\Ric}{Ric}
\DeclareMathOperator{\symrank}{symrank}

\DeclareMathOperator{\II}{I\! I}

\DeclareMathOperator{\trace}{trace}

\mathcode`l="8000
\begingroup
\makeatletter
\lccode`\~=`\l
\DeclareMathSymbol{\lsb@l}{\mathalpha}{letters}{`l}
\lowercase{\gdef~{\ifnum\the\mathgroup=\m@ne \ell \else \lsb@l \fi}}%
\endgroup

\newcommand*{\defeq}{\mathrel{\vcenter{\baselineskip0.5ex \lineskiplimit0pt
			\hbox{\scriptsize.}\hbox{\scriptsize.}}}%
	=}

\begin{document}

\title[Symmetry Rank and intermediate Ricci curvatures]{Local symmetry rank bound for \\ positive intermediate Ricci curvatures}
\author[Lawrence Mouill\'e]{Lawrence Mouill\'e}
\address{Department of Mathematics \\ Rice University \\ Houston, TX.}
\email{mouille@rice.edu}
\date{\today}
\blfootnote{\textup{2010} \textit{Mathematics Subject Classification}: Primary: 53C20; Secondary: 53B20, 58D19}
\blfootnote{Data sharing not applicable to this article as no datasets were generated or analyzed during the current study.}

\begin{abstract}
	We use a local argument to prove if an $r$-dimensional torus acts isometrically and effectively on a connected $n$-dimensional manifold which has positive $k^\mathrm{th}$-intermediate Ricci curvature at some point, then $r \leq \lfloor \frac{n+k}{2} \rfloor$.
	This symmetry rank bound generalizes those established by Grove and Searle for positive sectional curvature and Wilking for quasipositive curvature.
	As a consequence, we show that the symmetry rank bound in the Maximal Symmetry Rank Conjecture for manifolds of non-negative sectional curvature holds for those which also have positive intermediate Ricci curvature at some point.
	In the process of proving our symmetry rank bound, we also obtain an optimal dimensional restriction on isometric immersions of manifolds with non-positive intermediate Ricci curvature into manifolds with positive intermediate Ricci curvature, generalizing a result by Otsuki.
\end{abstract}

\maketitle

\section{Introduction}

	For manifolds with positive sectional curvature, the Grove symmetry program has been effective in inspiring innovative results. 
	The goal of this program, initiated by Karsten Grove in the 1990's, is to classify positively curved manifolds that have large isometry groups. 
	This initiative has led to advances in global Riemannian geometry, providing symmetry obstructions for positive curvature and motivating constructions of new examples of positive or non-negative curvature. 
	One prevalent measure of ``large isometry group'' is the symmetry rank:
	
	\begin{definition}
		The \textit{symmetry rank} of a Riemannian manifold $(M,g)$, which we denote by $\symrank(M,g)$, is the rank of its isometry group.
		In other words, $\symrank(M,g)$ the maximal dimension of a torus that can act isometrically and effectively on $(M,g)$. 
	\end{definition}
	
	A quintessential result from the Grove symmetry program is the ``Maximal Symmetry Rank Theorem'' established by Grove and Searle:
	
	\begin{theorem}[Grove \& Searle \cite{maxsymrank}]\label{thm:maxsymrank}
		If $(M^n,g)$ is a closed, connected, $n$-dimensional Riemannian manifold with positive sectional curvature, then $\symrank(M^n,g)\leq\left\lfloor\frac{n+1}{2}\right\rfloor$.
		Furthermore, if $\symrank(M^n,g)=\left\lfloor\frac{n+1}{2}\right\rfloor$, then $M^n$ is diffeomorphic to $S^n$, $\mathbb{R}\mathrm{P}^n$, $\mathbb{C}\mathrm{P}^{n/2}$, or a lens space.
	\end{theorem}
	
	Wilking generalized the symmetry rank bound from Theorem \ref{thm:maxsymrank} for manifolds with ``quasipositive curvature,'' i.e. non-negatively curved manifolds that contain a point at which all $2$-planes have positive sectional curvature.
	This symmetry rank does not appear in any of Wilking's publication; the author learned of this result through an unpublished manuscript prepared by Kerin. 
	Wilking's argument does not depend on the assumption of non-negative curvature, so we state his symmetry rank bound as follows:
	
	\begin{theorem}[Wilking \cite{wilkingunpublished}]\label{thm:quasipositive}
		Suppose $(M^n,g)$ is a connected, $n$-dimensional Riemannian manifold that contains a point at which all $2$-planes have positive sectional curvature. 
		Then $\symrank(M^n,g)\leq\lfloor\frac{n+1}{2}\rfloor$.
	\end{theorem}
	
	Readers can find Wilking's argument for Theorem \ref{thm:quasipositive} in Galaz-Garcia's article \cite{symrankquasiposlowdim}, in which Galaz-Garc\'ia establishes a diffeomorphism classification of $4$ and $5$-dimensional manifolds with quasipositive curvature and maximal symmetry rank.
	Notice that $M^n$ is not assumed to be closed in Theorem \ref{thm:quasipositive}, in contrast to Theorem \ref{thm:maxsymrank}.
	This is because Grove and Searle's argument is global in nature, depending on a result by Berger in \cite{BergerKilling} stating that any Killing field on a closed, even-dimensional manifold with globally positive sectional curvature must have a zero, while Wilking's argument is local in nature, depending only observations about the Gauss equation of an isometric immersion at a point where the extrinsic curvatures are positive.
	In this article, inspired by the Grove symmetry program and Theorems \ref{thm:maxsymrank} and \ref{thm:quasipositive}, we study the symmetry rank of manifolds with positive intermediate Ricci curvatures.
	
	\begin{definition}
		Given $n\geq 2$ and $k\in\{1,\dots,n-1\}$, an $n$-dimensional Riemannian manifold $(M^n,g)$ has \textit{positive $k^\mathrm{th}$-intermediate Ricci curvature at a point $p$} if the sum of sectional curvatures $\sum_{i=1}^k\sec(u,e_i)$ is positive for all orthonormal vectors $u,e_1,\dots,e_k$ at $p$.
		We will abbreviate this condition as ``$M$ has $\Ric_k|_p>0$.''
	\end{definition}
	
	For examples of manifolds with positive intermediate Ricci curvature and a survey of known structure results, see Section \ref{sec:context} below.
	Notice $\Ric_1|_p>0$ is equivalent to having positive sectional curvature at $p$, and $\Ric_{n-1}|_p>0$ is equivalent to having positive Ricci curvature at $p$.
	Furthermore, if $\Ric_k|_p>0$ for some $k$, then $\Ric_l|_p>0$ for all $l\geq k$.

	\subsection{Symmetry rank bound for positive intermediate Ricci curvatures}
		
		Our main result is the following symmetry rank bound, which generalizes the bounds from Theorems \ref{thm:maxsymrank} and \ref{thm:quasipositive}:
		
		\begin{Atheorem}\label{thm:global}
			If $(M^n,g)$ is a connected Riemannian $n$-manifold and contains a point $p$ at which $\Ric_k|_p>0$ for some $k\in\{1,\dots,n-1\}$, then
			\[
				\symrank(M^n,g)\leq\left\lfloor\frac{n+k}{2}\right\rfloor.
			\]
		\end{Atheorem}

		Our argument for Theorem \ref{thm:global} is inspired by Wilking's proof of Theorem \ref{thm:quasipositive}.
		Because ``$\Ric_1|_p>0$'' means that all $2$-planes at $p$ have positive sectional curvature, Theorem \ref{thm:quasipositive} is the $k=1$ case of Theorem \ref{thm:global}.
		Incidentally, similar to Theorem \ref{thm:quasipositive}, notice we do not assume $M^n$ is closed in Theorem \ref{thm:global}, in contrast to Theorem \ref{thm:maxsymrank}.
		
		The author recently sharpened the bound from Theorem \ref{thm:global} for closed manifolds which have $\Ric_k|_p>0$ at every point $p$ in \cite{TorusRic_k}.
		We mention one brief observation for this setting:
		Suppose $(M^n,g)$ is a closed, connected $n$-dimensional Riemannian manifold with $\Ric_k|_p > 0$ at every point $p$ for $k \geq n-2$.
		By Theorem \ref{thm:global}, $\symrank(M^n,g)\leq n-1$.
		The extremal case $\symrank(M^n,g)=n-1$ can be achieved in dimensions $n=2$ and $3$; e.g. the standard linear actions of $T^1$ on $S^2$ and $T^2$ on $S^3$, respectively.
		However, $\symrank(M^n,g)=n-1$ cannot be achieved for such manifolds in dimensions $n\geq 4$.
		This is because closed manifolds of dimension $n\geq 4$ which admit smooth effective torus actions of rank $n-1$ must have infinite fundamental group; see the work of Pak \cite{pak} and Parker \cite{parker}.
		Thus, by the Bonnet-Myers theorem, these manifolds cannot admit invariant metrics of positive Ricci curvature.
		Therefore, if $(M^n,g)$ is a closed, connected Riemannian manifold of dimension $n\geq4$ with $\Ric_k|_p > 0$ at every point $p$ for $k\geq n-2$, then $\symrank(M^n,g) \leq n-2$, which is a more restrictive bound than the one provided by Theorem \ref{thm:global}.
		In dimensions $n=4,5$, and $6$, products of spheres provide examples of closed Ricci-positive manifolds with $\symrank(M^n,g)=n-2$.
		Corro and Galaz-Garc\'ia show in \cite{riccisymrank} that in all dimensions $n \geq 5$, there exist closed $n$-manifolds with positive Ricci curvature and $\symrank(M^n,g)=n-4$. 
		Thus, we still have the following:
		
		\begin{question}
			Given a natural number $n\geq 7$, does there exist a closed, connected, $n$-dimensional Riemannian manifold $(M^n,g)$ with globally positive Ricci curvature such that $\symrank(M^n,g) = n-2$?
		\end{question}
	
	\subsection{Ramifications in non-negative curvature} 
		
		Theorem \ref{thm:global} has consequences related to the \textit{``Maximal Symmetry Rank Conjecture''} for non-negatively curved manifolds, which was formulated by Galaz-Garc\'ia and Searle in \cite{maxsymrankconj} and has since been sharpened by Escher and Searle in \cite{eschersearle}:
		
		\begin{conjecture}[Escher, Galaz-Garc\'ia, \& Searle \cite{eschersearle,maxsymrankconj}]\label{conj:maxsymrank}
			If $(M^n,g)$ is a closed, simply connected, $n$-dimensional Riemannian manifold with non-negative sectional curvature, then the following hold:
			\begin{enumerate}
				\item $\symrank(M^n,g)\leq\left\lfloor\frac{2n}{3}\right\rfloor$. \label{conjbound}
				
				\item If $\symrank(M^n,g)=\left\lfloor\frac{2n}{3}\right\rfloor$, then $M^n$ is equivariantly diffeomorphic to a product of spheres or a quotient thereof by a free linear action of a torus of rank less than or equal to $2n \;\mathrm{mod}\; 3$. 
				\label{diffeoclass}
			\end{enumerate}
		\end{conjecture}
		
		Galaz-Garc\'ia and Searle show in \cite{maxsymrankconj} that Conjecture \ref{conj:maxsymrank} holds for $4\leq n\leq 6$.
		Escher and Searle extend this result to dimensions $7\leq n\leq 9$ in \cite{eschersearle}.
		Furthermore, they show that Part (\ref{conjbound}) holds for $10\leq n\leq 12$, and they confirm Conjecture \ref{conj:maxsymrank} in all dimensions for torus actions that are \textit{``isotropy-maximal.''}
		With the assumption of non-negative curvature replaced with \textit{``rational ellipticity,''} Part (\ref{conjbound}) was established in all dimensions by Galaz-Garc\'ia, Kerin, and Radeschi in \cite{toruselliptic}. 
		Furthermore, they prove that Part (\ref{diffeoclass}) holds with ``equivariantly diffeomorphic'' replaced by ``rationally homotopy equivalent.'' 
		Their result is relevant to Conjecture \ref{conj:maxsymrank} because the Bott Conjecture claims that any non-negatively curved manifold is rationally elliptic.
		
		By Theorem \ref{thm:global}, we have the following immediate consequence:
		
		\begin{Acorollary}\label{cor:maxsymrankconj}
			Part (\ref{conjbound}) of Conjecture \ref{conj:maxsymrank} holds under the additional assumption that $M^n$ contains a point $p$ at which $\Ric_k|_p>0$ for $k=\lfloor n/3\rfloor$.
		\end{Acorollary}
		
		To justify Corollary \ref{cor:maxsymrankconj}, first notice that a manifold with non-negative sectional curvature can have $\Ric_k|_p>0$ at some point for $k\geq 2$ without having positive sectional curvature.
		For example, the Riemannian product $S^2\times S^2$ has $\Ric_3|_p > 0$ at every point $p$, but it does not have $\sec>0$; see also Example \ref{ex:product} below.
		Now, if a manifold $M^n$ has $\Ric_k>0$ at a point for $k\leq\frac{n}{3}$, then by Theorem \ref{thm:global}, $\symrank(M^n,g)\leq\lfloor\frac{2n}{3}\rfloor$.

	\subsection{Curvature restrictions on submanifolds}
		
		In the process of proving Theorem \ref{thm:global}, we produce a dimensional restriction on isometric immersions of submanifolds in terms of intermediate Ricci curvatures, which we will describe here.
		
		As mentioned before, Grove and Searle proved their symmetry rank bound in Theorem \ref{thm:maxsymrank} using a global argument, while Wilking established Theorem \ref{thm:quasipositive} using a local argument.
		Specifically, Wilking considered a principal orbit of an isometric torus action on a manifold at a point where the sectional curvatures are positive, studied the Gauss equation of its embedding, and obtained a dimensional restriction on the embedding.
		This restriction arising from Wilking's study of the Gauss equation is similar to a result, established by Otsuki in \cite{Otsuki}, which states if $N$ is non-positively curved, $N\hookrightarrow M$ is an isometric immersion, and $M$ is positively curved, then $\dim N \leq \lfloor\frac{\dim M + 1}{2}\rfloor$.
		
		In proving Theorem \ref{thm:global}, we generalize Otsuki's result to the following:
		
		\begin{Aproposition}\label{prop:immersion}
			Let $N$ be a Riemannian manifold that has $\Ric_k|_p\leq 0$ at a point $p$ for some $k\in\{1,\dots,\dim N-1\}$.
			Suppose $f:N\hookrightarrow M$ is an isometric immersion into a Riemannian manifold $M$ that has $\sum_{i=1}^k\sec(f_*(u),f_*(e_i)) > 0$ for all orthonormal vectors $u,e_1,\dots,e_k$ tangent to $N$ at $p$.
			Then 
			\[
				\dim N \leq \left\lfloor \frac{\dim M +k}{2}\right\rfloor.
			\]
		\end{Aproposition}
		
		By ``$\Ric_k|_p\leq 0$,'' we of course mean that $\sum_{i=1}^k\sec(u,e_i)\leq 0$ for all orthonormal $u,e_1,\dots,e_k\in T_pN$.
		For a more general version of Proposition \ref{prop:immersion}, see Corollary \ref{cor:Gauss} below.
		We will show that Proposition \ref{prop:immersion} is optimal in the following sense: 
		Given natural numbers $n \geq 2$, $k\in\{1,\dots,n-1\}$, and $m\geq n+1$ such that $n = \lfloor\frac{m+k}{2}\rfloor$, there exists an isometric immersion of Riemannian manifolds $f:N^n\hookrightarrow M^m$ and a point $p\in N^n$ such that $N^n$ has $\Ric_k|_p\leq 0$ and $M^m$ has $\sum_{i=1}^k\sec(f_*(u),f_*(e_i)) > 0$ for all orthonormal vectors $u,e_1,\dots,e_k\in T_pN$.
		For a more precise statement, see Proposition \ref{prop:models} below.

\subsection{Structure of this article}
In Section \ref{sec:context}, we present examples of manifolds with positive intermediate Ricci curvature and survey some results from the literature. 
In Section \ref{sec:proof}, we prove Proposition \ref{prop:immersion} and use it to establish Theorem \ref{thm:global}. 
In Section \ref{sec:example}, we construct isometric immersions which demonstrate that Proposition \ref{prop:immersion} is optimal.

\subsection{Acknowledgments}

The results in this article are part of my doctoral thesis at the University of California, Riverside, and many were completed while I was supported by the university's Dissertation Year Program Award. 
I thank my thesis advisor, Fred Wilhelm, for his insights and valuable feedback during the preparation of this article. 
I also thank Priyanka Rajan, Catherine Searle, Karsten Grove, and Guofang Wei for helpful discussions concerning this project, and an anonymous reviewer for their detailed feedback on a previous version of this article.

\section{Examples \& context}\label{sec:context}

	Many have worked to find new examples of manifolds with positive sectional curvature, non-negative sectional curvature, positive Ricci curvature, etc.
	However, only a few have published examples of manifolds with positive intermediate Ricci curvature.
	We mention a few elementary sources of examples here:
	
	\begin{example}\label{ex:product}
		Let $M^m$ and $N^n$ be manifolds of dimension $m$ and $n$, respectively, and suppose each has positive sectional curvature at points $p$ and $q$, respectively.
		Now consider the point $(p,q)$ in their Riemannian product $M^m\times N^n$.
		Because the Levi-Civita connection of $M^m\times N^n$ splits as the sum of the connections of the factors, it is a simple exercise to show that $M^m\times N^n$ has $\Ric_k|_{(p,q)} > 0$ only for $k \geq \max\{m,n\}+1$.
		For example, given $n\geq 2$, the $2n$-dimensional Riemannian product of spheres $S^n\times S^n$ has $\Ric_{n+1}|_{(p,q)}>0$ at every point $(p,q)$.
	\end{example}
	
	\begin{example}\label{ex:symmetric}
		For the classification of compact, irreducible, symmetric spaces according to the minimal values of $k$ for which each has $\Ric_k|_p>0$ at every point $p$, see \cite{AQZ} or \cite{DGM}.
		For example, a compact, irreducible, symmetric space has positive sectional curvature at every point if and only if it is of rank $1$.
		Among those of rank $\geq 2$, none have $\Ric_2|_p>0$ at every point $p$, and the only ones with $\Ric_3|_p>0$ for every point $p$ are the $5$-dimensional space $\mathsf{SU}_3/\mathsf{SO}_3$, the $6$-dimensional space $\mathsf{SO}_5/(\mathsf{SO}_2\times\mathsf{SO}_3)$, and the $8$-dimensional space $\mathsf{G}_2/\mathsf{SO}_4$.
	\end{example}
	
	\begin{example}\label{ex:S^3xS^3}
		The product of spheres $S^3\times S^3$ is diffeomorphic to the quotient of $S^3\times S^3\times S^3$ by the diagonal $S^3$-action.
		When $S^3\times S^3\times S^3$ is equipped with the standard Riemannian product metric, the diagonal $S^3$-action is by isometries, and hence the quotient $S^3\times S^3$ inherits a Riemannian metric such that the quotient map $S^3\times S^3\times S^3 \to S^3 \times S^3$ is a Riemannian submersion.
		This metric on $S^3 \times S^3$ has $\Ric_2|_p>0$ at every point $p$, which is a stronger curvature condition than what holds for the Riemannian product metric discussed in Example \ref{ex:product}.
		For a generalization of this observation to products of compact semisimple Lie groups, see Theorem E in \cite{DGM}.
	\end{example}
	
	We close this section by referencing several existing structure results for manifolds with globally positive intermediate Ricci curvature.
	First, note if a compact manifold has $\Ric_k|_p>0$ at every point $p$, then it has positive Ricci curvature, and by the Bonnet-Myers theorem, must have finite fundamental group.
	Many results from the setting of positive or non-negative sectional curvature have been generalized for intermediate Ricci curvature.
	These include generalizations of the Synge theorem and Weinstein fixed point theorem by Wilhelm in \cite{wilhelminterm}, the Heintze-Karcher inequality by Chahine in \cite{yousef}, and the Gromoll-Meyer theorem and Cheeger-Gromoll soul theorem for open manifolds by Shen in \cite{shenkth}.
	In addition, many comparison results have been proven by Guijarro and Wilhelm in their series of papers \cite{softconnprinc,focalradius,restrictions}.
	For connectedness principles of submanifolds, see Remark 2.4 in Wilking's article \cite{connprinc}, Theorem C in Fang, Mendon{\c{c}}a, and Rong's article \cite{FangMendoncaRong}, or Theorem E the author's article \cite{TorusRic_k}.
	In \cite{TorusRic_k}, the author uses a connectedness principle to classify closed manifolds that have $\Ric_2|_p>0$ at every point $p$ and maximal symmetry rank, generalizing Theorem \ref{thm:maxsymrank}.
	Gumaer and Wilhelm produce a restriction on dual foliations in \cite{GumaerWilhelm}, generalizing a result by Wilking from \cite{transversejacobi}.
	Finally, for a characterization of positive intermediate Ricci curvature through optimal transport, see Section 6 in Ketterer and Mondino's article \cite{KettererMondino}.

\section{Isometric immersions \& symmetry rank bound}
\label{sec:proof}
	
	In this section, we prove Proposition \ref{prop:immersion} and Theorem \ref{thm:global} using arguments inspired by Wilking's proof of Theorem \ref{thm:quasipositive}.
	The idea is to consider an isometric immersion of a torus into a manifold of $\Ric_k > 0$ (e.g. consider a principal orbit of an isometric torus action), and study the Gauss equation for the curvature.
	Because the torus has intrinsic curvature zero, the extrinsic curvature condition $\Ric_k > 0$ produces a restriction on the Gaussian curvature of the immersion, which in turn produces a restriction on the dimension of the torus.
	
	\subsection{Isometric immersions}
		
		First, we discuss isometric immersions and intermediate Ricci curvatures in order to prove Proposition \ref{prop:immersion}.
		Let $f:N \hookrightarrow M$ be an isometric immersion of Riemannian manifolds, and let $\II$ denote the second fundamental form for the immersion.
		In other words, for each $p\in N$, $\II:T_pN \times T_pN \to T_p^\perp N$ is the symmetric bilinear form defined by $\II(u,v)=(\nabla_uV)^\perp\in T_p^\perp N \subset T_pM$, where $V$ is any extension of $v$ to a vector field.
		Then the Gauss equation asserts that the intrinsic sectional curvature of $N$, $\sec_N$, is related to sectional curvature of $M$, $\sec_M$, as follows: 
		For all orthonormal $u,v\in T_pN$,
		\begin{equation}\label{eq:Gauss}
			\sec_M(f_*(u),f_*(v)) = \sec_N(u,v) + |\II(u,v)|^2-\langle \II(u,u),\II(v,v)\rangle.
		\end{equation}
		Studying the algebraic properties of Gaussian curvatures $|\II(u,v)|^2-\langle \II(u,u),\II(v,v)\rangle$ leads to restrictions on the dimensions of $N$ and $M$.
		For example, Otsuki proved in \cite{Otsuki} that if a symmetric bilinear form $B:V\times V\to W$ satisfies $|B(u,v)|^2-\langle B(u,u),B(v,v)\rangle > 0$ for all orthonormal $u,v$, then $\dim W \geq \dim V-1$.
		In particular, by applying this result to the second fundamental form $\II:T_pN \times T_pN \to T_p^\perp N$, if the Gaussian curvature of an isometric immersion $f:N\hookrightarrow M$ is always positive, then $\dim M \geq 2\dim N - 1$.
		Here, we generalize Otsuki's result as follows:
		
		\begin{lemma}\label{lem:Otsuki}
			Let $V$ and $W$ be real finite-dimensional vector spaces endowed with inner products, let $B:V\times V\to W$ be a symmetric bilinear form, and let $k\in\{1,\dots,\dim V - 1\}$. 
			If 
			\begin{equation}\label{eq:bilinearform}
				\sum_{i=1}^k\left[|B(u,e_i)|^2-\langle B(u,u),B(e_i,e_i)\rangle\right]>0
			\end{equation}
			for all orthonormal sets of vectors $u,e_1,\dots,e_k \in V$, then $\dim W\geq \dim V-k$.
		\end{lemma}
		
		\begin{proof}
			For the sake of contradiction, assume $\dim W\leq \dim V-k-1$.
			We will show there exist orthonormal vectors $u,e_1,\dots,e_k \in V$ such that Inequality \eqref{eq:bilinearform} does not hold.
			First, let $S$ denote the unit sphere in $V$, and define the function $f:S\to\mathbb{R}$ by 
			\[
			f(v)=|B(v,v)|^2. 
			\]
			Notice $f$ is a continuous function with compact domain, so $f$ has a minimal value. 
			Choose a vector $u \in S$ at which $f$ achieves its minimum.
			Let $\mathrm{span}\{u\}^\perp$ denote the orthogonal complement of $u$ in $V$, and define the operator $B_{u}:\mathrm{span}\{u\}^\perp \to W$ by $B_{u}(v) = B(u,v)$.
			By assumption, $\dim \mathrm{span}\{u\}^\perp = \dim V - 1 \geq  \dim W + k$.
			Thus, by the rank-nullity theorem,
			\[
				\dim \ker B_{u} \geq \dim \mathrm{span}\{u\}^\perp -\dim W \geq k.
			\]
			In particular, we can choose orthonormal vectors $e_1,\dots,e_k\in\mathrm{span}\{u\}^\perp$ such that $B_u(e_i) = B(u,e_i)=0$ for each $i$.
			So it suffices to show that $\langle B(u,u),B(e_i,e_i)\rangle \geq 0$ for each $i$.
			Now, for $i\in\{1,\dots,k\}$, define the family of vectors $u_i:\mathbb{R}\to S$ by
			\[
			u_i(\theta)\defeq \cos(\theta) u + \sin(\theta) e_i.
			\]
			Notice that $u_i(0) = u$.
			Because $f$ is minimal at $u$, for each $i$, we have
			\begin{align*}
				0 &\leq \tfrac{d^2}{d\theta^2}f(u_i(\theta))\big|_{\theta=0}\\
				&=4(\langle B(u,u),B(e_i,e_i)\rangle -|B(u,u)|^2).
			\end{align*}
			Thus, it follows that $\langle B(u,u),B(e_i,e_i)\rangle \geq |B(u,u)|^2 \geq 0$.
			Therefore, because $B(u,e_i)=0$ for each $i$, we have that 
			\[
			\sum_{i=1}^k\left[|B(u,e_i)|^2-\langle B(u,u),B(e_i,e_i)\rangle\right]\leq 0.\qedhere
			\]
		\end{proof}
		
		Related to the $k=1$ case of Lemma \ref{lem:Otsuki}, Otsuki also proved the following algebraic lemma, which was conjectured by Chern and Kuiper in \cite{ChernKuiper}:
		If there exists $\lambda \geq 0$ such that a symmetric bilinear form $B:V\times V\to W$ satisfies $|B(u,u)| > \lambda$ and $|B(u,v)|^2-\langle B(u,u),B(v,v)\rangle \geq \lambda^2$ for all orthonormal $u,v$, then $\dim W \geq \dim V$.
		One consequence of this lemma is if $M^n$ is a compact, non-positively curved manifold that admits an isometric immersion into a Euclidean space $\mathbb{R}^{n+d}$, then $d\geq n$. 
		For a survey of this and related results, see Chapter 10 Section 3 in \cite{DoCarmo} or Section 3.1 in \cite{Dajczer}.
		Rovenski\u{i} proved an analogue of this algebraic lemma in \cite{Rovenskii} and used it to establish restriction results on isometric immersions of manifolds whose domains have upper bounds on intermediate Ricci curvatures.
		For example, Rovenski\u{i} showed if $M^n$ is a compact manifold that has $\Ric_k|_p\leq 0$ at every point $p$ and $M^n$ admits an isometric immersion into a Euclidean space $\mathbb{R}^{n+d}$, then $d\geq n - k$. 
		
		We will now use the Gauss equation and Lemma \ref{lem:Otsuki} to establish the following:
		
		\begin{corollary}\label{cor:Gauss}
			Let $f:N \hookrightarrow M$ be an isometric immersion of Riemannian manifolds and $k\in\{1,\dots,\dim N-1\}$. 
			Suppose there exists a point $p\in N$ and a subspace $\mathcal{S} \subseteq T_pN$ such that $\dim \mathcal{S} \geq k+1$ and 
			\[
				\sum_{i=1}^k \sec_N(u,e_i) < \sum_{i=1}^k \sec_M(f_*(u),f_*(e_i))
			\]
			for all orthonormal vectors $u,e_1,\dots,e_k$ in $\mathcal{S}$.
			Then $\dim \mathcal{S} \leq \dim M -\dim N + k$.
		\end{corollary}
		
		\begin{proof}
			Given any orthonormal vectors $u,e_1,\dots,e_k$ in $\mathcal{S}$, by Equation \eqref{eq:Gauss},
			\[
			0< \sum_{i=1}^k [\sec_M(f_*(u),f_*(e_i))-\sec_N(u,e_i)]=\sum_{i=1}^k[|\II(u,e_i)|^2-\langle \II(u,u),\II(e_i,e_i)\rangle].
			\]
			Let $\II|_\mathcal{S}:\mathcal{S} \times \mathcal{S} \to T_p^\perp N$ denote the second fundamental form restricted to the subspace $\mathcal{S}$.
			Then Corollary \ref{lem:Otsuki} applies to the symmetric bilinear form $\II|_\mathcal{S}$, implying that $\dim M - \dim N \geq \dim \mathcal{S} - k$.
		\end{proof}
		
		Now Proposition \ref{prop:immersion} easily follows from Corollary \ref{cor:Gauss}:
		
		\begin{proof}[Proof of Proposition \ref{prop:immersion}]
			Let $f:N \hookrightarrow M$ be an isometric immersion of Riemannian manifolds, let $k\in\{1,\dots,\dim N-1\}$, and suppose there exists a point $p\in N$ such that $N$ has $\Ric_k|_p\leq 0$ and $M$ has $\sum_{i=1}^k\sec(f_*(u),f_*(e_i)) > 0$ for all orthonormal vectors $u,e_1,\dots,e_k\in T_pN$.
			Thus, for all orthonormal vectors $u,e_1,\dots,e_k\in T_pN$,
			\[
				\sum_{i=1}^k \sec_N(u,e_i) < \sum_{i=1}^k \sec_M(f_*(u),f_*(e_i)).
			\]
			Then by Corollary \ref{cor:Gauss} with $\mathcal{S} = T_pN$, it follows that $\dim N \leq \lfloor \frac{\dim M + k}{2} \rfloor$.
		\end{proof}
		
	\subsection{Symmetry rank bound}
		
		We will now prove Theorem \ref{thm:global}.
		As mentioned at the beginning of this section, we will apply Proposition \ref{prop:immersion} to the principal orbits of an isometric torus action.
		First, we justify that the manifolds under consideration indeed must have positive intermediate Ricci curvature at a point on some principal orbit:
		
		\begin{lemma}\label{lem:principal}
			Suppose a torus $T$ acts smoothly on a connected Riemannian manifold $M$.
			If $M$ has $\Ric_k|_p>0$ at some point $p$ on a singular orbit of the $T$-action, then there exists a point $\tilde{p}$ that lies on a principal orbit such that $M$ has $\Ric_k|_{\tilde{p}}>0$.
		\end{lemma}
		
		\begin{proof}
			First, we describe how $\Ric_k$ can be viewed as a real-valued function defined on partial flags within the tangent bundle of $M^n$.
			Let $\mathrm{Fl}(1,k+1;TM)$ denote the {partial flag bundle} consisting of signature-$(1,k+1)$ flags tangent to $M$. 
			In other words, elements of $\mathrm{Fl}(1,k+1;TM)$ are pairs $(\mathcal{V}^1,\mathcal{V}^{k+1})$ such that $\mathcal{V}^1$ is a $1$-dimensional subspace of a $(k+1)$-dimensional subspace $\mathcal{V}^{k+1}$ of a tangent space $T_pM$ for some $p\in M$.
			Given such a flag $(\mathcal{V}^1,\mathcal{V}^{k+1})$ and a choice of unit vector $u\in\mathcal{V}^1$, let $J_u|_{\mathcal{V}^{k+1}}$ denote the {Jacobi operator} restricted to $\mathcal{V}^{k+1}$.
			In other words, letting $R$ denote the type-$(1,3)$ curvature tensor for $(M,g)$ and letting $P_{\mathcal{V}^{k+1}}:T_pM \to \mathcal{V}^{k+1}$ denote the orthogonal projection onto $\mathcal{V}^{k+1}$, define $J_u|_{\mathcal{V}^{k+1}}:\mathcal{V}^{k+1}\to \mathcal{V}^{k+1}$ by
			\[
			J_u|_{\mathcal{V}^{k+1}}(x) \defeq P_{\mathcal{V}^{k+1}}(R(x,u)u).
			\]
			Then we can define the function $\Ric_k\colon \mathrm{Fl}(1,k+1;TM)\to\mathbb{R}$ by
			\[
			\Ric_k(\mathcal{V}^1,\mathcal{V}^{k+1})\defeq\trace\big(J_u|_{\mathcal{V}^{k+1}}\big)=\sum_{i=1}^k\sec(u,e_i),
			\]
			where $u$ is any unit vector in the line $\mathcal{V}^1$, and $e_1,\dots,e_k$ is any orthonormal basis for the orthogonal complement of $\mathcal{V}^1$ in $\mathcal{V}^{k+1}$. 
			This map is well-defined because the value of $\Ric_k(\mathcal{V}^1,\mathcal{V}^{k+1})$ is independent of the choice of unit vector $u\in\mathcal{V}^1$ and orthonormal vectors $e_1,\dots,e_k\in\mathcal{V}^{k+1}\cap\mathcal{V}^1$.
			
			Now suppose that $M$ has $\Ric_k|_p>0$ at some point $p$ that lies on a singular orbit of the $T$-action on $M$.
			Because the $T$-action on $M$ is smooth and $M$ is connected, the set of principal orbits is dense in $M$; see, for example, Chapter IV Section 3 in \cite{Bredon}.
			Hence, the set of partial flags within $\mathrm{Fl}(1,k+1;TM)$ that are tangent to principal orbits is dense in $\mathrm{Fl}(1,k+1;TM)$.
			Notice, because $\Ric_k\colon \mathrm{Fl}(1,k+1;TM)\to\mathbb{R}$ is defined using the curvature tensor and orthogonal projections, it is continuous.
			Because $\Ric_k$ is continuous and $\Ric_k|_p>0$, there exists an open trivializing neighborhood for the bundle $\mathrm{Fl}(1,k+1;TM)$ around $p$ on which $\Ric_k > 0$.
			Therefore, because the set of partial flags tangent to principal orbits is dense in $\mathrm{Fl}(1,k+1;TM)$, there exists a point $\tilde{p}$ near $p$ such that $\tilde{p}$ lies on a principal orbit and $\Ric_k|_{\tilde{p}}>0$. 
		\end{proof}
		
		We are now ready to apply Proposition \ref{prop:immersion} and Lemma \ref{lem:principal} to prove Theorem \ref{thm:global}:
		
		\begin{proof}[Proof of Theorem \ref{thm:global}]
			Let $(M^n,g)$ be a connected, $n$-dimensional Riemannian manifold, suppose $\Ric_k|_p>0$ at a point $p\in M^n$ for some $k\in\{1,\dots,n-1\}$, and assume a torus $T^r$ acts isometrically and effectively on $M^n$. 
			We must show that $r  \leq \lfloor\frac{n+k}{2}\rfloor$.
			
			First, we address the case when $k=n-1$.
			Because $\lfloor\frac{n+(n-1)}{2}\rfloor = n-1$, it suffices in this case to show that a torus $T^n$ cannot act isometrically and effectively $M^n$ while $M^n$ has positive Ricci curvature at some point $p$.
			Suppose $T^n$ does act isometrically and effectively on $M^n$.
			Because the action is effective, the principal isotropy groups are trivial, and hence the principal orbits are $n$-dimensional.
			Thus, the tangent spaces of $M^n$ at every point of every principal orbit is spanned by commuting Killing fields.
			It then follows that $M^n$ has zero sectional curvature on all principal orbits.
			But because the set of principal orbits are dense in $M^n$, it follows that $M^n$ has zero sectional curvature everywhere.
			Therefore, if $T^n$ acts isometrically and effectively on $M^n$, then $M^n$ cannot have positive Ricci curvature at any point, thus establishing the case $k=n-1$.
			
			Now assume $k\leq n-2$.
			Let $N$ denote the $T^r$-orbit through $p$ in $M^n$.
			By Lemma \ref{lem:principal}, we may assume that $N$ is a principal orbit.
			Again, because the $T^r$-action is effective, $N$ is $r$-dimensional.
			Notice if $k \geq r$, then because $r = \dim N \leq \dim M = n$, it follows that $r \leq \lfloor\frac{n+k}{2}\rfloor$.
			So suppose $k \leq r-1$.
			Then, by applying Proposition \ref{prop:immersion} to the immersion $N \hookrightarrow M^n$, we have that $r \leq \lfloor\frac{n+k}{2}\rfloor$.
		\end{proof}

\section{Isometric immersions of maximal Euclidean rank}
\label{sec:example}
	
	In Proposition \ref{prop:immersion} we established if a manifold $N$ with $\Ric_k\leq 0$ is isometrically immersed in a manifold $M$ with $\Ric_k>0$, then $\dim N \leq \lfloor \frac{\dim M + k}{2}\rfloor$.
	In this section, we prove that this inequality is optimal.
	Specifically, we establish the following:
	
	\begin{proposition}\label{prop:models}
		Suppose $m,n,k$ are natural numbers such that $2 \leq n$, $n+1 \leq m$, $1\leq k\leq n-1$, and $n = \lfloor \frac{m+k}{2}\rfloor$.
		Then there exists a metric $g_{m,n,k}$ on $\mathbb{R}^m$ such that the flat Euclidean space $\mathbb{R}^n$ admits an isometric immersion $f:\mathbb{R}^n\hookrightarrow (\mathbb{R}^m,g_{m,n,k})$, and for every set of orthonormal vectors $u,e_1,\dots,e_k$ tangent to $\mathbb{R}^n$ at the origin, $\mathbb{R}^m$ has $\sum_{i=1}^k\sec(f_*(u),f_*(e_i))>0$.
	\end{proposition}
	
	Assume $m$, $n$, and $k$ are natural numbers that satisfy the hypotheses of Proposition \ref{prop:models}.
	The immersion $f:\mathbb{R}^n\hookrightarrow \mathbb{R}^m$ will be given by the inclusion $(y_1,\dots,y_n)\mapsto(0,\dots,0,y_1,\dots,y_n)$, and the metrics $g_{m,n,k}$ will be warped product metrics on $\mathbb{R}^{m-n}\times \mathbb{R}^n$.
	We will make it so that the coordinate vector fields for the $\mathbb{R}^n$ factor are Killing fields, meaning that $f$ will be an isometric immersion, and hence the image of $\mathbb{R}^n$ under $f$ will have zero intrinsic curvature.
	In order to ensure $\mathbb{R}^m$ has $\sum_{i=1}^k\sec(f_*(u),f_*(e_i))>0$ for all orthonormal vectors $u,e_1,\dots,e_k$ tangent to $\mathbb{R}^n$ at the origin, we will use the following computational simplification: 
	
	\begin{lemma}\label{lem:compsimp}
		Let $M^m$ be an $m$-dimensional Riemannian manifold, let $p$ be a point in $M^m$, and let $\mathcal{N}\subseteq T_pM$ be an $n$-dimensional subspace with orthonormal basis $Y_1,\dots,Y_n$ for some $n\in\{1,\dots,m-1\}$.
		Assume $k\in\{1,\dots,n-1\}$ and let $\mu$ and $\nu$ be non-negative real numbers.
		Suppose the following hold:
		\begin{enumerate}[(A)]
			\item $R(Y_i,Y_j)Y_l=0$ when the indices $i,j,l$ are mutually distinct.
			\label{a}
			
			\item $\sec(Y_i,Y_j)\in\{-\nu,\mu\}$ for all $i\neq j$.
			\label{b}

			\item For each $i$, there exist at most $k-1$ indices $j$ such that $j\neq i$ and $\sec(Y_i,Y_j)=-\nu$.
			\label{c}
			
			\item $\mu-(k-1)\nu>0$.
			\label{d}
		\end{enumerate}
		Then $M^m$ has $\sum_{i=1}^k\sec(u,e_i)>0$ for all orthonormal vectors $u,e_1,\dots,e_k\in\mathcal{N}$.
	\end{lemma}
	
	\begin{proof}
		First, we establish a lower bound on Ricci curvatures.
		Let $\Ric|_{\mathcal{N}}:\mathcal{N}\to\mathcal{N}$ denote the Ricci $(1,1)$-tensor restricted to $\mathcal{N}$ and composed with projection onto $\mathcal{N}$. 
		Then by definition, $\langle \Ric|_{\mathcal{N}}(u),u\rangle=\sum_{i=1}^{n-1}\sec(u,e_i)$ for any choice of orthonormal vectors $e_1,\dots,e_{n-1}\in\mathcal{N}$ that are orthogonal to $u$. By Property \eqref{a}, we get that $\Ric|_{\mathcal{N}}$ is diagonalized by $\{Y_i\}_{i=1}^n$; see, for example, \cite[Proposition 4.1.3]{petersen}. Then by Properties \eqref{b} and \eqref{c}, it follows that for any orthonormal basis $u,e_1,\dots,e_{n-1}$ of $\mathcal{N}$, 
		\begin{align*}
			\sum_{i=1}^{n-1}\sec(u,e_i) =\langle\Ric|_{\mathcal{N}}(u),u\rangle &\geq \min_{i=1,\dots,n}\{\langle\Ric|_{\mathcal{N}}(Y_i),Y_i\rangle\} \\
			&= \min_{i=1,\dots,n}\Big\{\sum_{j\neq i}\sec(Y_i,Y_j)\Big\} \\
			&\geq (n-k)\mu-(k-1)\nu.
		\end{align*}
		Applying Property \eqref{d}, it follows that
		\begin{equation}
			\sum_{i=1}^{n-1}\sec(u,e_i)>(n-k-1)\mu.
			\label{riccibound}
		\end{equation}
		
		Now, we establish an upper bound on sectional curvature.
		Define the operator $\mathfrak{R}|_{\mathcal{N}}:\Lambda^2\mathcal{N}\to \Lambda^2\mathcal{N}$ by $\langle\mathfrak{R}|_{\mathcal{N}}\left(\textstyle\sum_i x_i\wedge z_i\right),\textstyle\sum_j v_j\wedge w_j\rangle=\sum_{i,j} R(x_i,z_i,w_j,v_j)$; i.e. the curvature operator restricted to
		the subspace $\mathcal{N}$. 
		Then by Property \eqref{a}, $\mathfrak{R}|_{\mathcal{N}}$ is diagonalized by $\{Y_i\wedge Y_j\}_{i<j}$; see, for example, \cite[Proposition 4.1.2]{petersen}. 
		Thus by Property \eqref{b}, for every pair of orthonormal vectors $u,v\in\mathcal{N}$, we have 
		\begin{equation}
			\sec(u,v) = \langle\mathfrak{R}|_{\mathcal{N}}(u\wedge v),u\wedge v\rangle \leq \max_{i<j}\{\sec(Y_i\wedge Y_j)\} = \mu.
			\label{secbound}
		\end{equation}	
		
		Finally, consider any $k+1$ orthonormal vectors $u,e_1,\dots,e_k\in\mathcal{N}$.
		To this list of vectors, append any $n-k-1$ orthonormal vectors $e_{k+1},\dots,e_{n-1}\in\mathcal{N}$ such that $\{u,e_1,\dots,e_k,e_{k+1},\dots,e_{n-1}\}$ is an orthonormal basis for $\mathcal{N}$. 
		Then by Inequalities \eqref{riccibound} and \eqref{secbound}, we have
		\begin{align*}
			\sum_{i=1}^k\sec(u,e_i) &= \sum_{i=1}^{n-1}\sec(u,e_i)-\sum_{i=k+1}^{n-1}\sec(u,e_i)\\
			&> \left[(n-k-1)\mu\right]-\left[(n-k-1)\mu\right]=0.\qedhere
		\end{align*}
	\end{proof}
	
	Now we use Lemma \ref{lem:compsimp} to construct the metrics from Proposition \ref{prop:models}:
	
	\begin{proof}[Proof of Proposition \ref{prop:models}]
		Suppose $m,n,k$ are fixed natural numbers such that $2 \leq n$, $n+1 \leq m$, $1\leq k\leq n-1$, and $n = \lfloor \frac{m+k}{2}\rfloor$.
		We will construct a warped product metric $g_{m,n,k}$ on $\mathbb{R}^m = \mathbb{R}^{m-n}\times \mathbb{R}^n$ such that the inclusion $\mathbb{R}^n\hookrightarrow \mathbb{R}^{m-n}\times \mathbb{R}^n$ into the second factor is an isometric immersion of the flat Euclidean space $\mathbb{R}^n$, and for every set of orthonormal vectors $u,e_1,\dots,e_k$ tangent to the $\mathbb{R}^n$-factor of $\mathbb{R}^{m-n}\times \mathbb{R}^n$ at the origin, $\mathbb{R}^m$ has $\sum_{i=1}^k\sec(u,e_i)>0$.
		
		Consider coordinates $(x_1,\dots,x_{m-n},y_1,\dots,y_n)$ on $\mathbb{R}^{m-n}\times \mathbb{R}^n$. 
		Given positive smooth functions $\phi_i:\mathbb{R}^{m-n}\to\mathbb{R}$, we define the metric $g_{m,n,k}$ on $\mathbb{R}^{m-n}\times\mathbb{R}^n$ by
		\[
			g_{m,n,k}={dx_1}^2+\dots+{dx_{m-n}}^2+{\phi_1}^2{dy_1}^2+\dots+{\phi_n}^2{dy_n}^2.
		\]
		Notice that the the coordinate vector fields $\frac{\partial}{\partial y_1},\dots,\frac{\partial}{\partial y_n}$ are commuting Killing fields under the metric $g_{m,n,k}$.
		Thus, the inclusion of $\mathbb{R}^n$ into the second factor of $\mathbb{R}^{m-n}\times \mathbb{R}^n$ is in fact an isometric immersion of the flat Euclidean space $\mathbb{R}^n$.
		Define orthonormal vectors $Y_1,\dots,Y_n$ tangent to $\mathbb{R}^m$ at the origin by
		\[
			Y_i\defeq\tfrac{1}{\phi_{i}(0)}\tfrac{\partial}{\partial y_i}\big|_0.
		\]
		We will choose the functions $\phi_1,\dots,\phi_n$ such that $\left(\mathbb{R}^m,g_{m,n,k}\right)$ and $Y_1,\dots,Y_n$ together satisfy the hypotheses of Lemma \ref{lem:compsimp}. 
		
		Let $\mathcal{N}=\mathrm{span}\{Y_i\}_{i=1}^d=T_0\mathbb{R}^n$, and denote the orthogonal complement by $\mathcal{N}^\perp=T_0\mathbb{R}^{m-n}$. Let $\II:\mathcal{N}\times\mathcal{N}\to \mathcal{N}^\perp$ denote the second fundamental form for the submanifold $\{(0,\dots,0)\}\times\mathbb{R}^n$. Then $\II(Y_i,Y_j)=0$ for $i\neq j$, and combining this with the fact that $\frac{\partial}{\partial y_i}$ are commuting Killing fields, we get that Property \eqref{a} of Lemma \ref{lem:compsimp} is satisfied. 
		For Properties \eqref{b}, \eqref{c}, and \eqref{d} to be satisfied, it suffices to choose the functions $\phi_1,\dots,\phi_n$ so that \begin{enumerate}[(a)]
			\item $-\nu \defeq \sec(Y_i,Y_j) = \sec(Y_i,Y_l) < 0$ for all distinct $i,j,l\in\{1,\dots,k\}$,  \label{aa}
			
			\item \hspace{6pt}$\mu \defeq \sec(Y_i,Y_\alpha) = \sec(Y_\beta,Y_\gamma) > 0$ for all $i\in\{1,\dots,k\}$ and all $\alpha,\beta,\gamma \in \{k+1,\dots,n\}$ such that $\beta \neq \gamma$,\label{bb}
			
			\item $\mu-(k-1)\nu > 0$.\label{cc}
		\end{enumerate}
		Because $\II(Y_i,Y_i)=-\frac{\nabla \phi_i}{\phi_i(0)}$ and the submanifold $\{(0,\dots,0)\}\times\mathbb{R}^n$ has zero intrinsic curvature, it follows from the Gauss equation that the extrinsic curvatures are given by 
		\begin{equation}\label{eq:sectional}
			\sec(Y_i,Y_j)=-\frac{\langle\nabla \phi_i,\nabla \phi_j\rangle}{\phi_i(0)\phi_j(0)}, \quad \text{for }i,j\in\{1,\dots,n\},\quad i\neq j.
		\end{equation}
		For the sake of simplicity, we will choose $\phi_1,\dots,\phi_n$ such that
		\begin{enumerate}
			\item $a \defeq \phi_1(0) = \dots = \phi_k(0) > 0$,\label{11}
			
			\item $b \defeq \phi_{k+1}(0) = \dots = \phi_n(0) > 0$,\label{22}
			
			\item $\nabla \phi_1 = \dots = \nabla \phi_k$,\label{33}
			
			\item $\langle \nabla \phi_\alpha , \nabla \phi_\beta \rangle = \langle \nabla \phi_\gamma , \nabla \phi_\delta \rangle < 0$ \\ 
			for all $\alpha,\beta,\gamma,\delta\in\{k+1,\dots,n\}$ such that $\alpha\neq \beta$ and $\gamma\neq \delta$, \\ 
			(assuming $k+1 \neq n$)\label{44}
			
			\item $\langle \nabla \phi_i , \nabla \phi_\alpha \rangle = \langle \nabla \phi_j , \nabla \phi_\beta \rangle < 0$ \\
			for all $i,j\in\{1,\dots,k\}$ and $\alpha,\beta\in\{k+1,\dots,n\}$.\label{55}
		\end{enumerate}
		It follows from \eqref{eq:sectional} and \eqref{11}--\eqref{55} that the negative sectional curvatures among $\sec(Y_i,Y_j)$  are those for which $i,j\in\{1,\dots,k\}$, $i\neq j$ (assuming $k\geq 2$; otherwise, there are no negative curvatures).
		These negative curvatures have the value $-\nu = -\frac{|\nabla \phi_1|^2}{a^2}$.
		All other sectional curvatures are positive, and they either have value $-\frac{\langle \nabla \phi_1, \nabla \phi_n \rangle}{ab}$ or $-\frac{\langle \nabla \phi_{k+1} , \nabla \phi_n \rangle}{b^2}$ (assuming $k+1\neq n$).
		So in order to satisfy \eqref{bb}, meaning that the positive curvatures among $\sec(K_i,K_j)$ all have the same value, we will require that
		\begin{enumerate}\setcounter{enumi}{5}
			\item $\frac{\langle \nabla \phi_1, \nabla \phi_n \rangle}{ab} = \frac{\langle \nabla \phi_{k+1} , \nabla \phi_n \rangle}{b^2}$ (assuming $k+1 \neq n$).\label{66}
		\end{enumerate}
		Finally, to satisfy \eqref{cc}, we must also require 
		\begin{enumerate}\setcounter{enumi}{6}
			\item $ - \frac{\langle \nabla \phi_1, \nabla \phi_n \rangle}{ab} - (k-1)\frac{| \nabla \phi_1|^2}{a^2}> 0$.\label{77}
		\end{enumerate}
		In summary, our goal is to define $\phi_1,\dots,\phi_n$ such that Conditions \eqref{11}--\eqref{77} are satisfied, and Proposition \ref{prop:models} then follows from Lemma \ref{lem:compsimp}.
		Recall that $\nabla \phi_1,\dots,\nabla \phi_n\in \mathcal{N}^\perp = T_0\mathbb{R}^{m-n}$.
		Our choices of $\phi_1,\dots,\phi_n$ will be rather different depending on whether $m-n = 1$ or $m-n\geq 2$.
		Thus, we now break the construction into two cases: $m-n = 1$ and $m - n \geq 2$. 
		
		\textbf{Case $\mathbf{m-n=1}$:} 
			If $m-n = 1$, then $\mathcal{N}^\perp = T_0\mathbb{R}^{m-n}$ is $1$-dimensional, so $\phi_1,\dots,\phi_n$ are functions of just one variable. 
			Also, because $n=\lfloor\frac{m+k}{2}\rfloor$, it follows that $k+1=n$, so we may ignore Conditions \eqref{44} and \eqref{66}.
			Define $U$ to be a choice of unit vector in $\mathcal{N}^\perp$.
			To satisfy Conditions \eqref{11}, \eqref{22}, \eqref{33}, and \eqref{55}, for some constants $a,b>0$, we will choose $\phi_1,\dots,\phi_n:\mathbb{R}\to\mathbb{R}$ to be any smooth functions such that
			\begin{itemize}
				\item $\phi_1(0)=\dots=\phi_k(0)=a$,
				
				\item $\nabla\phi_1=\dots=\nabla\phi_k=U$,
				
				\item $\phi_n(0)=b$, and
				
				\item $\nabla\phi_n=-U$.
			\end{itemize}
			Finally, to satisfy Condition \eqref{77}, we must choose $a,b>0$ so that 
			\[
				- \frac{\langle \nabla \phi_1, \nabla \phi_n \rangle}{ab} - (k-1)\frac{| \nabla \phi_1|^2}{a^2} = \frac{1}{ab} - \frac{k-1}{a^2} > 0
			\]
			Therefore, choosing any values for $a$ and $b$ such that $a>(k-1)b$, the $m-n=1$ case of Proposition \ref{prop:models} follows from Lemma \ref{lem:compsimp}.
			
		\textbf{Case $\mathbf{m-n \geq 2}$:} 
			Suppose now that $m-n\geq 2$, meaning $\mathcal{N}^\perp = T_0\mathbb{R}^{m-n}$ is at least $2$-dimensional.
			Choose any unit vector $U\in\mathcal{N}^\perp$.
			To satisfy Conditions \eqref{11}, \eqref{22}, and \eqref{33}, we will choose $\phi_1,\dots,\phi_n$ such that
			\begin{enumerate}[(i)]
				\item $\phi_1(0)=\dots=\phi_k(0)=a$ for some $a>0$ to be determined, \label{i}
				
				\item $\phi_{k+1}(0)=\dots=\phi_n(0)=b$ for some $b>0$ to be determined,
				
				\item $\nabla\phi_1=\dots=\nabla\phi_k=U$.
			\end{enumerate}
			To satisfy Conditions \eqref{44} and \eqref{55}, we will choose $\nabla\phi_{k+1},\dots,\nabla\phi_n$ to be vectors $V_{k+1}(\theta), \dots, V_n(\theta)$ such that the angle between $U$ and $V_i(\theta)$ is independent of the choice of index $i$, and the angle between $V_i(\theta)$ and $V_j(\theta)$ is independent of the choice of distinct indices $i$ and $j$.
			To that end, consider the unit sphere $S^{m-n-1}\subset \mathcal{N}^\perp$, and let $S^{m-n-2}\subset S^{m-n-1}$ denote the equatorial sphere that lies in the subspace orthogonal to $U$. 
			Consider any inscribed regular $(m-n-1)$-simplex in $S^{m-n-2}$, and define $V_{k+1}(0),\dots,V_n(0)$ to be the vertices of that simplex. 
			This is possible because $n-k-1$, the number of vectors $V_{k+1}(0),\dots,V_n(0)$, is equal to $m-n$, the number of vertices of the $(m-n-1)$-simplex.
			The inner products of these vectors is given by
			\begin{equation}\label{eq:innerprodV}
				\langle V_i(0),V_j(0)\rangle=-\frac{1}{m-n-1}
			\end{equation}
			for all distinct $i,j\in\{k+1,\dots,n\}$.
			Now, given $\theta\in[0,\frac{\pi}{2}]$, we define
			\begin{equation}\label{eq:defV}
				V_i(\theta)\defeq\cos(\theta)V_i(0)-\sin(\theta)U.
			\end{equation}
			See Figure \ref{figure} for an illustration of $U,V_{k+1}(\theta),\dots,V_n(\theta)$ in $S^{m-n-1}$.
			
			\begin{figure}[H]
				\centering
				\includegraphics[width=\textwidth]{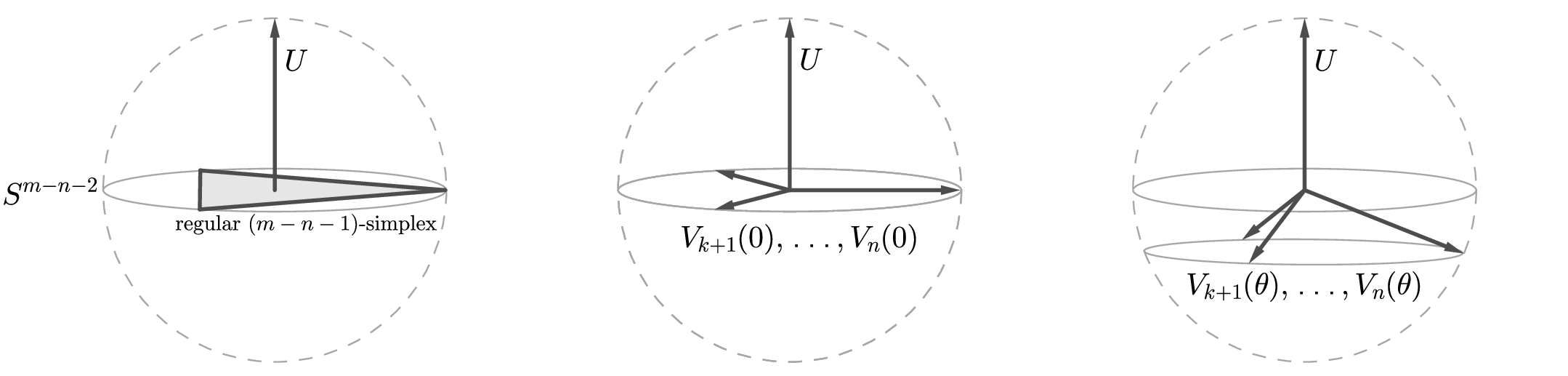}
				\caption{The vectors $U,V_{k+1}(\theta),\dots,V_n(\theta)$ in $S^{m-n-1}$}
				\label{figure}
			\end{figure}
			
			It follows from Equations \eqref{eq:innerprodV} and \eqref{eq:defV} that for all distinct $i$ and $j$, 
			\begin{equation}\label{eq:<Vi,Vj>}
				\langle V_i(\theta),V_j(\theta)\rangle=\sin^2\theta-\big(\tfrac{1}{m-n-1}\big)\cos^2\theta.
			\end{equation}
			Notice this function of $\theta$ is strictly increasing on $[0,\frac{\pi}{2}]$, takes the value $-\frac{1}{m-n-1}$ at $\theta=0$, and takes the value $0$ at 
			\[
				\theta_0=\arctan\left(\sqrt{\tfrac{1}{m-n-1}}\right).
			\]
			Hence, $\langle V_i(\theta),V_j(\theta)\rangle$ is negative for $\theta\in[0,\theta_0)$.
			Thus, choosing $\phi_{k+1},\dots,\phi_n$ as follows, Conditions \eqref{44} and \eqref{55} will be satisfied:
			\begin{enumerate}[(i)]\setcounter{enumi}{3}
				\item $\nabla \phi_i = V_i(\theta)$ as defined in \eqref{eq:defV}, for all $i\in\{k+1,\dots,n\}$, for some $\theta\in(0,\theta_0)$ to be determined,\label{iv}
			\end{enumerate} 
			We will now demonstrate that for functions $\phi_1,\dots,\phi_n$ satisfying \eqref{i}--\eqref{iv}, there exists values of $a,b>0$ and $\theta\in(0,\theta_0)$ such that Conditions \eqref{66} and \eqref{77} are satisfied.
			In particular, we must show $a,b,\theta$ can be chosen such that
			\begin{equation}\label{eq:=}
				\frac{\langle U, V_n(\theta) \rangle}{ab} = \frac{\langle V_{k+1}(\theta) , V_n(\theta) \rangle}{b^2},
			\end{equation}
			\begin{equation}\label{eq:>}
				- \frac{\langle U, V_n(\theta) \rangle}{ab} - (k-1)\frac{|U|^2}{a^2}> 0.
			\end{equation}
			To satisfy Equation \eqref{eq:=}, we will define $a>0$ so that
			\begin{enumerate}[(i)]\setcounter{enumi}{4}
				\item $a=\frac{b\langle U, V_n(\theta) \rangle}{\langle V_{k+1}(\theta) , V_n(\theta) \rangle}$ for some $b>0$ and $\theta\in(0,\theta_0)$ to be determined.
			\end{enumerate}
			Notice that this value of $a$ is indeed positive because the inner products $\langle U, V_n(\theta) \rangle$ and $\langle V_{k+1}(\theta) , V_n(\theta) \rangle$ are both negative for $\theta\in(0,\theta_0)$.
			Using this value of $a$, Inequality \eqref{eq:>} then reduces to 
			\begin{equation}\label{eq:>reduced}
				\langle U,V_n(\theta)\rangle^2+(k-1)|U|^2\langle V_{k+1}(\theta),V_n(\theta)\rangle > 0.
			\end{equation}
			In particular, notice this inequality does not depend on $b$.
			Thus, we can choose $b$ so that
			\begin{enumerate}[(i)]\setcounter{enumi}{5}
				\item $b$ is any positive real number.
			\end{enumerate}
			Finally, because $|U| = 1$, it follows from Equations \eqref{eq:defV} and \eqref{eq:<Vi,Vj>} that Inequality \eqref{eq:>reduced} is equivalent to 
			\[
				k \sin^2\theta - \tfrac{k-1}{m-n-1} \cos^2\theta > 0.
			\]
			Notice that the quantity on the left side of this inequality is continuous in the variable $\theta$ and takes the value $\frac{1}{m-n}$ when $\theta = \theta_0 = \arctan(\sqrt{1/(m-n-1)})$.
			Thus, there exist a value of $\theta\in(0,\theta_0)$ such that the quantity on the left is positive.
			Hence, we can indeed choose $\theta$ so that
			\begin{enumerate}[(i)]\setcounter{enumi}{6}
				\item $\theta\in(0,\theta_0)$ such that $k \sin^2\theta - \tfrac{k-1}{m-n-1} \cos^2\theta > 0$.\label{vii}
			\end{enumerate}
			Therefore, defining functions $\phi_1,\dots,\phi_n:\mathbb{R}^{m-n}\to\mathbb{R}$ to be any smooth functions such that \eqref{i}--\eqref{vii} hold, Proposition \ref{prop:models} follows from Lemma \ref{lem:compsimp}.
	\end{proof}

\bibliographystyle{abbrv}
\bibliography{symmetry_rank_biblio}

\end{document}